\documentclass[final,3p,times]{elsarticle}
\usepackage{amsfonts}

 \usepackage{graphics}
 \usepackage{graphicx}
 \usepackage{epsfig}

\usepackage{amssymb}
\usepackage{amsmath,latexsym}
\usepackage{mathrsfs}
\usepackage{bm}
\usepackage{amssymb}
\usepackage{hyperref}
\usepackage{hyperref}
\newtheorem{thm}{Theorem}
\newtheorem{lem}{Lemma}
\newtheorem{pro}{Proposition}
\newdefinition{de}{Definition}
\newdefinition{rmk}{Remark}
\newproof{pf}{Proof}
\newproof{pot}{Proof of the existence of Theorem 3.2}
\newtheorem{hypo}{Hypothesis}

\journal{Bernoulli}

\begin{document}

\begin{frontmatter}



\title{Chaotic and Predictable Representations for
Multidimensional L\'{e}vy Processes\tnoteref{t1}}
\tnotetext[t1]{This research was supported by the National Basic
Research Program of China (973 Program) (Program No.2007 CB814903)
and the National Natural Science Foundation of China (Program
No.70671069).}

\author[sjtu]{Jianzhong Lin\corref{cor1}}

\ead{jzlin@sjtu.edu.cn} \cortext[cor1]{Corresponding author.}

\address[sjtu]{Department of Mathematics, Shanghai Jiaotong University, Shanghai 200240, China}

\begin{abstract}
For a general Multidimensional L\'{e}vy process (satisfying some
moment conditions), we introduce the Multidimensional power jump
processes and the related Multidimensional Teugels martingales.
Furthermore, we orthogonalize the Multidimensional Teugels
martingales by applying Gram-Schmidt process. We give a chaotic
representation for every square integral random variable in terms of
these orthogonalized Multidimensional Teugels martingales. The
predictable representation with respect to the same set of
Multidimensional orthogonalized martingales of square integrable
random variables and of square integrable martingales is an easy
consequence of the chaotic representation.
\end{abstract}


\end{frontmatter}


\section{Introduction}\label{sec1}

The \rm{chaotic representation property(CRP)} has been studied by
Emery(1989) for normal martingales, that is, for martingales $X$
such that $<X,X>_t=ct$, for some constant $c>0$. This property says
that any square integrable random variable measurable with respect
to $X$ can be expressed as an orthogonal sum of multiple stochastic
integrals with respect to $X$. It is known (see for example
Dellacherie et al., 1992, p. 207 and Dermoune, 1990), that the only
normal martingales $X$, with the CRP, or even the weaker
\rm{predictable representation property (PRP)}, which are also
L\'{e}vy processes are the Brownian motion and the compensated
Poisson process. David Nualart and Wim Schoutens (2000) study the
chaotic representation property for one-dimensional L\'{e}vy
process, in terms of a suitable orthogonal sequence of martingales
where these martingales are obtained as the orthogonalization of the
compensated power jump processes of the L\'{e}vy process.
Furthermore, Nualart and Schoutens (2001) used their martingale
representation result to establish the existence and uniqueness of
solutions for BSDE's driven by a L\'{e}vy process of the kind
considered in Nualart and Schoutens (2000). In addition, Corcuera,
Nualart and Schoutens (2005) applied this martingale representation
result to study the completion of a L\'{e}vy market by power-jump
assets.

In the past twenty years, there is already a growing interest for
multidimensional L\'{e}vy Processes. Some concepts and basic
properties about multidimensional L\'{e}vy Processes were summarized
in Sato (1999). Applications of multidimensional L\'{e}vy Processes
to analyzing biomolecular (DNA and protein) data and one-server
light traffic queues were explored by Dembo, Karlin and Zeittouni
(1994). A small deviations property of multidimensional L\'{e}vy
Processes were discussed by Simon (2003). In finance research,
practically all financial applications require a multivariate model
with dependence between components: examples are basket option
pricing, portfolio optimization, simulation of risk scenarios for
portfolios. In most of these applications, jumps in the price
process must be taken into account. Cont and Tankov (2004)
systematically investigated these problems in multidimensional
L\'{e}vy market. In addition, the optimal portfolios in
multidimensional L\'{e}vy market is discussed by Emmer and
Kl\"{u}ppelberg (2004). Some simulation approaches for multivariate
L\'{e}vy processes are also investigated in Cohen and Rosi\'{n}ski
(2007). L\'{e}vy copulas was also suggested by Kallsen and Tankov
(2006) in order to characterize the dependence among components of
multidimensional L\'{e}vy Processes. The SDEs driven by
infinite-dimensional L\'{e}vy processes was investigated by
Meyer-Brandis and Proske (2010).

The chaotic representation property is important for the research of
L\'{e}vy process, and multidimensional L\'{e}vy processes obtain
some applications in bioscience and finance, so it is significant to
extend the result in univariate set-up obtained by Nualart and
Schouten(2000) to the cases of multidimensional L\'{e}vy processes.
In this paper, following the research line of the paper in Nualart
and Schouten(2000), we study the chaotic representation property for
Multi-dimensional L\'{e}vy processes, in terms of a suitable
orthogonal sequence of Multidimensional martingales, assuming that
the L\'{e}vy measure has a finite Laplace transform outside the
origin. These Multidimensional martingales are obtained as the
orthogonalization of the Multidimensional compensated power jump
processes of our Multidimensional L\'{e}vy process. In Section 2, we
introduce these Multidimensional compensated power jump processes
and we transform them into a multivariate orthogonal sequence.
Section 3 is devoted to prove the chaos representation property from
which a predictable representation is deduced. Finally, in Section
4, we discuss some particular examples.

\section{Preliminary}\label{sec2}

A $\mathbb{R}^n$-valued stochastic process
$X=\{X(t)=(X_1(t),X_2(t),\cdots,X_n(t))',t\geq 0\}$ defined in
complete probability space $(\Omega,\mathscr{F},\mathbb{P})$ is
called \textsl{L\'{e}vy process} if $X$ has stationary and
independent increments and $X(0)=\bm{0}$. A L\'{e}vy process
possesses a c\`{a}dl\`{a}g modification (Protter,1990, Theorem
30,p.21) and we will always assume that we are using this
c\`{a}dl\`{a}g version. If we let
$\mathscr{F}_t=\mathscr{G}_t\vee\mathscr{N}$, where
$\mathscr{G}_t=\sigma\{X(s),0\leq s\leq t\}$ is the natural
filtration of $X$, and $\mathscr{N}$ are the $\mathbb{P}-$null sets
of $\mathscr{F}$, then $\{\mathscr{F}_t,t\geq 0\}$ is a right
continuous family of $\sigma-$fields (Protter,1990,Theorem 31,p.22).
We assume that $\mathscr{F}$ is generated by $X$. For an up-to-date
and comprehensive account of L\'{e}vy processes we refer the reader
to Bertoin (1996) and Sato (1999).

Let $X$ be a L\'{e}vy process and denote by
\begin{eqnarray}
X(t-)=\lim\limits_{s\rightarrow t,s<t}X(s),\quad t>0 ,\nonumber
\end{eqnarray}
the left limit process and by $\triangle X(t)=X(t)-X(t-)$ the jump
size at time $t$. It is known that the law of $X(t)$ is
\textsl{infinitely divisible} with characteristic function of the
form
\begin{eqnarray}
E\left[exp(i\bm{\theta}\cdot
X(t))\right]=\left(\phi(\bm{\theta})\right)^t,\quad
\bm{\theta}=(\theta_1,\theta_2,\cdots,\theta_n)\in
\mathbb{R}^n\nonumber
\end{eqnarray}
where $\phi(\bm{\theta})$ is the characteristic function of
$\bm{X}(1)$. The function $\psi(\bm{\theta})=log\phi(\bm{\theta})$
is called the \textsl{characteristic exponent} and it satisfies the
following famous L\'{e}vy-Khintchine formula (Bertoin, 1996):
\begin{eqnarray}
\psi(\bm{\theta})=-\frac{1}{2}\bm{\theta}\cdot
\Sigma\bm{\theta}+\textrm{i}\bm{a}\cdot
\bm{\theta}+\int_{\mathbb{R}^n}\left(exp(i\bm{\theta}\cdot
\bm{x})-1-\textrm{i}\bm{\theta}\cdot \bm{x} 1_{|\bm{x}|\leq
1}\right)\nu(d\bm{x}).\nonumber
\end{eqnarray}
where $\bm{a},\bm{x}\in \mathbb{R}^n$, $\Sigma$ is a symmetric
nonnegative-definite $n\times n$ matrix, and $\nu$ is a measure on
$\mathbb{R}^n\backslash\{o\}$ with $\int(\|\bm{x}\|^2\wedge
1)\nu(d\bm{x})<\infty$. The measure $\nu$ is called the
\textsl{L\'{e}vy measure} of $X$.

Throughout this paper, we will use the standard multi-index
notation. We denote by $\mathbb{N}_0$ the set of nonnegative
integers. A multi-index is usually denoted by $\bm{p}$,
$\bm{p}=(p_1,p_2,\cdots,p_n)\in\mathbb{N}_0^n$. Whenever $\bm{p}$
appears with subscript or superscript, it means a multi-index. In
this spirit, for example, for $\bm{x}=(x_1,\cdots,x_n)$, a monomial
in variables $x_1,\cdots,x_n$ is denoted by
$\bm{x}^{\bm{p}}=x_1^{p_1}\cdots x_n^{p_n}$.  In addition, we also
define $\bm{p}!=p_1!\cdots p_n!$ and $|\bm{p}|=p_1+\cdots+p_n$; and
if $\bm{p}$, $\bm{q}\in\mathbb{N}_0^n$, then we define
$\delta_{\bm{p},\bm{q}}=\delta_{\bm{p}_1,\bm{q}_1}\cdots\delta_{\bm{p}_n,\bm{q}_n}$.

\begin{hypo}
We will suppose in the remaining of the paper that the L\'{e}vy
measure satisfies for some $\varepsilon>0$, and $\lambda>0$,
\begin{eqnarray}
\int_{|\bm{x}|\geq \epsilon}exp(\lambda
\|\bm{x}\|)\nu(d\bm{x})<\infty. \nonumber
\end{eqnarray}
\end{hypo}
This implies that
\begin{eqnarray}
\int \bm{x}^{\bm{p}}\nu(d\bm{x})<\infty. \quad |\bm{p}|\geq 2
\end{eqnarray}
and that the characteristic function $E\left[exp(i\bm{\theta}\cdot
X(t))\right]$ is analytic in a neighborhood of origin $\bm{o}$. As a
consequence, $X(t)$ has moments of all orders and the polynomials
are dense in $L^2(\mathbb{R}^n,\mathbb{P}\circ X(t)^{-1})$ for all
$t>0$.

Professor Nualart proposed author to use the following
transformations of $X$ which will play an important role in our
analysis. We introduce power jump monomial processes of the form
\begin{eqnarray}
X(t)^{(p_1,\cdots,p_n)}\stackrel{\rm{def}}{=}\sum\limits_{0<s\leq
t}(\triangle X_1(s))^{p_1}\cdots(\triangle X_n(s))^{p_n},\nonumber
\end{eqnarray}
The number $|\bm{p}|$ is called the total degree of $X(t)^{\bm{p}}$.
Furthermore define
\begin{eqnarray}
Y(t)^{(p_1,\cdots,p_n)}\stackrel{\rm{def}}{=}X(t)^{(p_1,\cdots,p_n)}-\mathbb{E}[X(t)^{(p_1,\cdots,p_n)}]
=X(t)^{(p_1,\cdots,p_n)}-m_{\bm{p}}t ,\nonumber
\end{eqnarray}
the compensated power jump process of multi-index
$\bm{p}=(p_1,p_2,\cdots,p_n)$. Under hypothesis 1,
$Y(t)^{(p_1,\cdots,p_n)}$ is a normal martingale, since for an
integrable L\'{e}vy process $Z$, the process $\{Z_t-E[Z_t], t\geq
0\}$ is a martingale. We call $Y(t)^{(p_1,\cdots,p_n)}$ the
\textsl{Teugels martingale monomial} of multi-index
$(p_1,\cdots,p_n)$.

\begin{rmk}
In the case of a Poisson process, all power jump processes will be
the same, and equal to the original \textsl{Poisson process}. In the
case of a \textsl{Brownian motion}, all power jump processes of
order strictly greater than one will be equal to zero.
\end{rmk}

In the following we will introduce some concepts and basic
properties of martingale polynomial. These new concepts and
properties are totaly similar to those of polynomial in $n$ real
variables (cf. Dunkl and Xu (2001)).

A \textsl{martingale polynomial} $P$ in $n$ L\'{e}vy variables
$X=(X_1,X_2,\cdots,X_n)$ is a linear combination of martingale
monomials,
\begin{eqnarray}
P(X)=\sum\limits_{|\bm{p}|\geq 1}c_{\bm{p}}Y^{\bm{p}},\nonumber
\end{eqnarray}
where the coefficients $c_{\bm{p}}$ are in the real numbers
$\mathbb{R}$. The degree of a martingale polynomial is defined as
the highest total degree of its martingale monomials. We shall use
the abbreviation $\Pi^n$ to denote the collection of all martingale
polynomials in $X$. We also denote the space of martingale
polynomials of degree at most $d$ by $\Pi_d^n$. A martingale
polynomial is called \textsl{homogeneous} if all the monomials
appearing in it have the same total degree. Denote the space of
homogeneous polynomials of degree $d\in\mathbb{N}$ in $n$ variables
by $\mathcal{P}_d^n$; that is
\begin{eqnarray}
\mathcal{P}_d^n=\left\{P:P(X)=\sum\limits_{|\bm{p}|=d}c_{\bm{p}}Y^{\bm{p}}\right\}.\nonumber
\end{eqnarray}
Every polynomial in $\Pi^n$ can be written as a linear combination
of homogeneous martingale polynomials; for $P\in \Pi_d^n$,
\begin{eqnarray}
P(X)=\sum\limits_{k=1}^d\sum\limits_{|\bm{p}|=k}c_{\bm{p}}Y^{\bm{p}}.\nonumber
\end{eqnarray}
Denote by $r_d^n$ the dimension of $\mathcal{P}_d^n$ and it is well
known that
\begin{eqnarray}
r_d^n=dim\mathcal{P}_d^n=\left(\begin{array}{c}
                           d+n-1 \\
                           d
                         \end{array}\right)
&and&dim\Pi_d^n=\left(\begin{array}{c}
                           d+n \\
                           d
                         \end{array}\right)-1.\nonumber
\end{eqnarray}

We denote by $\mathscr{N}^2$ the space of one dimensional square
integrable martingales $M$ such that
$\sup\limits_t\mathbb{E}(M(t)^2)<\infty$, and $M(0)=0$ a.s. Notice
that if $M\in\mathscr{N}^2$, then
$\lim\limits_{t\rightarrow\infty}\mathbb{E}(M(t)^2)=\mathbb{E}(M(\infty)^2)<\infty$,
and $M(t)=\mathbb{E}[M(\infty)|\mathscr{F}_t]$. Thus, each $M\in
\mathscr{N}^2$ can be identified with its terminal value
$M(\infty)$. As in Protter(2005, p.181), we say that two martingales
$M,N\in\mathscr{N}^2$ are strongly orthogonal and we denote this by
$M\times N$, if and only if the product $MN$ is a uniformly
integrable martingale. As noted in Protter (2005,p181), one can
prove that $M\times N$ if only if $[M,N]$ is a uniformly integrable
martingale. We say that two random vectors $X,Y\in
L^2(\Omega,\mathscr{F})$ are weakly orthogonal, $X\perp Y$,if
$E[XY]=0$. Clearly, strong orthogonality implies weak orthogonality.

In the theory about orthogonal polynomials of several variables, we
can apply the Gram-Schmidt process to the monomials with respect to
the usual inner product to produce a sequence of orthogonal
polynomials of several variables. Some details about the technique
and theory of \textsl{orthogonal polynomials of several variables}
refer to Dunkl and Xu (2001). In this paper, we shall apply the
standard Gram-Schmidt process with the graded lexicographical order
to generate a biorthogonal basis
$\{H^{\bm{p}},\bm{p}\in\mathbb{N}^n\}$, such that each
$H^{\bm{p}}(|\bm{p}|=d)$ is a linear combination of the
$Y^{\bm{q}}\in\Pi_d^n$, with $|\bm{q}|\leq |\bm{p}|$ and the leading
coefficient equal to $1$. We set
\begin{eqnarray}
H^{\bm{p}}&=&Y^{\bm{p}}+\sum\limits_{\bm{q}\prec\bm{p},|\bm{q}|=|\bm{p}|}c_{\bm{q}}Y^{\bm{q}}+\sum\limits_{k=1}^{|\bm{p}|-1}\sum\limits_{|\bm{q}|=k}c_{\bm{q}}Y^{\bm{q}},\nonumber
\end{eqnarray}
where $\bm{p}=\{p_1,\cdots,p_n\}$, $\bm{q}=\{q_1,\cdots,q_n\}$ and
$\prec$ represent the relation of graded lexicographical order
between two multi-indexes.

We have that
\begin{eqnarray}
\begin{array}{rl}
[H^{\bm{p}},Y^{\bm{q}}](t)=&\sum\limits_{0<s\leq t}(\triangle
X_1(s))^{p_1+q_1}\cdots(\triangle
X_n(s))^{p_n+q_n}\\
&+\sum\limits_{0<s\leq t}\sum\limits_{1\leq
|\tilde{\bm{p}}|<|\bm{p}|}c_{\tilde{\bm{p}}+\bm{q}}(\triangle
X_1(s))^{\tilde{p}_1+q_1}\cdots(\triangle X_n(s))^{\tilde{p}_n+q_n}
+\sum\limits_{i=1}^n\sum\limits_{j=1}^nc_{\bm{e}_i+\bm{q}}\sigma_{ij}tI_{\{\bm{q}=\bm{e}_j\}}.\nonumber
\end{array}
\end{eqnarray}
Let
$m_{\bm{p}+\bm{q}}=\int\prod\limits_{i=1}^nx_i^{p_i+q_i}\nu(d\bm{x})$,
then we have that
\begin{eqnarray}
\mathbb{E}[H^{\bm{p}},Y^{\bm{q}}](t)&=&t\left(m_{\bm{p}+\bm{q}}
+\sum\limits_{1\leq
|\tilde{\bm{p}}|<|\bm{p}|}c_{\tilde{\bm{p}}+\bm{q}}m_{\tilde{\bm{p}}+\bm{q}}
+\sum\limits_{i=1}^n\sum\limits_{j=1}^nc_{\bm{e}_i+\bm{q}}\sigma_{ij}tI_{\{\bm{q}=\bm{e}_j\}}\right).\nonumber
\end{eqnarray}
where $\bm{e}_i$ denotes the n-dimensional unit vector with $i$th
component equal to one. In conclusion, we have that,
$[H^{\bm{p}},Y^{\bm{q}}]$ is a martingale if only if we have that
$E[H^{\bm{p}},Y^{\bm{q}}](1)=0$.

Consider two spaces: The first space $S_1$ is defined as follows
\begin{eqnarray}
S_1&=&\left\{\sum\limits_{k=1}^d\sum\limits_{|\bm{p}|=k}c_{k}(p_1,\cdots,p_n)x_1^{p_1}\cdots
x_n^{p_n}+c_0+\sum\limits_{(i_1,\cdots,i_n)\in\{0,-1\}^n,|\bm{i}|\geq
-(n-1)}c_{-1}(i_1,\cdots,i_n)x_1^{i_1}\cdots
x_n^{i_n}; \right.\nonumber\\
&&\left. d\in\{1,2,\cdots\},
c_{j}(p_1,\cdots,p_n)\in\mathbb{R},j=-1,0,\cdots,d; x_i\neq 0,
i=1,2,\cdots,n;\bm{i}=(i_1,\cdots,i_n)\right\}\nonumber
\end{eqnarray}
which is endowed with the scalar product $<\cdot,\cdot>_1$ given by
\begin{eqnarray}
<P(\bm{x}),Q(\bm{x})>_1&=&\int_{-\infty}^{+\infty}\cdots\int_{-\infty}^{+\infty}P(\bm{x})Q(\bm{x})
\prod\limits_{i=1}^nx_i^2\nu(d\bm{x})\nonumber\\
&&+\sum\limits_{i=1}^n\sum\limits_{j=1}^nc_1c_2\sigma_{ij}I_{\{P(\bm{x})=c_1\bm{x}^{\bm{e}_i-\bm{1}},
Q(\bm{x})=c_2\bm{x}^{\bm{e}_j-\bm{1}}\}},\qquad
\bm{1}=(1,1,\cdots,1).\nonumber
\end{eqnarray}
Note that
\begin{eqnarray}
&&<x_1^{p_1-1}\cdots x_n^{p_n-1},x_1^{q_1-1}\cdots x_n^{q_n-1}>_1\nonumber\\
&=&m_{\bm{p+q}}+\sum\limits_{i=1}^n\sum\limits_{j=1}^n\sigma_{ij}I_{\{\bm{p}=\bm{e}_i,\bm{q}=\bm{e}_j\}},\quad
|\bm{p}|\geq 1,\quad |\bm{q}|\geq 1 .\nonumber
\end{eqnarray}

Thus we can construct the other space $S_{2}$ which is the space of
all linear transformations of the Teugels martingale monomials of
the multivariate L\'{e}vy process, i.e.
\begin{eqnarray}
S_{2}&=\left\{\sum\limits_{p_1+\cdots+p_n=d}c_{d}(p_1,\cdots,p_n)Y(t)^{(p_1,\cdots,p_n)}
+\sum\limits_{p_1+\cdots+p_n=d-1}c_{d-1}(p_1,\cdots,p_n)Y(t)^{(p_1,\cdots,p_n)}\right.\nonumber\\
&\left.+\cdots+\sum\limits_{p_1+\cdots+p_n=1}c_{1}(p_1,\cdots,p_n)Y(t)^{(p_1,\cdots,p_n)},\quad
d\geq 1,\quad i=1,2,\cdots,n .\right\}.\nonumber
\end{eqnarray}
We endow this space with the scalar product $<\cdot,\cdot>_2$, given
by
\begin{eqnarray}
<Y^{(p_1,\cdots,p_n)},Y^{(q_1,\cdots,q_n)}>_2&=&E([Y^{(p_1,\cdots,p_n)},Y^{(q_1,\cdots,q_n)}](1))\nonumber\\
&=&m_{\bm{p+q}}+\sum\limits_{i=1}^n\sum\limits_{j=1}^nc_{\bm{e}_i+\bm{q}}\sigma_{ij}tI_{\{\bm{q}=\bm{e}_j\}},
\nonumber\\
&&|\bm{p}|\geq 1,\quad |\bm{q}|\geq 1 .\nonumber
\end{eqnarray}

Because $\sum\limits_{0<s\leq t}(\triangle
X_1(s))^{p_1}\cdots(\triangle
X_n(s))^{p_n}\equiv\sum\limits_{0<s\leq t}(\triangle
X_1(s))^{p_1}\cdots(\triangle X_n(s))^{p_n}I_{\cap_{i=1}^n(\triangle
X_i(s)\neq 0)}$, one clearly sees that

$x_1^{p_1-1}x_2^{p_2-1}\cdots x_n^{p_n-1}\leftrightarrow
Y^{(p_1,\cdots,p_n)}$ with $|\bm{p}|\geq 1$ and is an isometry
between $S_{1}$ and $S_{2}$. An orthogonalization of
$\{x_1^{-1}x_2^{-1}\cdots x_{n-1}^{-1},x_1^{-1}x_3^{-1}\cdots
x_n^{-1},\cdots,x_n^{-1},1,x_1,\cdots,x_n,x_1^2,x_1x_2,\cdots,x_n^2,\cdots\}$in
$S_{1}$ gives an orthogonalization of

$\{Y^{(1,0,\cdots,0)},\cdots,Y^{(0,\cdots,0,1)},Y^{(2,0,\cdots,0)},Y^{(1,1,0,\cdots,0)},\cdots,Y^{(0,\cdots,0,2)},\cdots\}$.

It is well known that \textsl{orthogonal polynomials of several
variables are not unique} (cf. Dunkl and Xu (2001)). In the
remaining of the paper, $\{H^{\bm{p}},\bm{p}\in\mathbb{N}^n\}$ is a
set of pairwise strongly orthogonal martingales given by the
previous orthogonalization of $\{Y^{\bm{p}},\bm{p}\}$. It is also
worth to emphasis that all deduction procedures and results are the
same once the orthogonal martingales are determinatively given.

\section{Representation properties}
\label{sec3}
\subsection{Representation of a power of a L\'{e}vy process}
For notation simplicity, here and hereafter we set
\begin{eqnarray}
X_i^{(p_i)}(t)=\sum_{0<s\leq t}(\triangle X_i(s))^{p_i},\quad
p_i\geq 2,\quad i=1,2,\cdots,n.\nonumber
\end{eqnarray}
and for convenience we put $X_i^{(1)}(t)=X_i(t)$. $p_i,
q_i(i=1,2,\cdots,n)$ are all some nonnegative integers. Note that
not necessarily $X_i(t)=\sum\limits_{0<s\leq t}\triangle X_i(s)$
holds; it is only true in the bounded variation case with
$\Sigma=O$. If $\Sigma=O$, clearly $[X_i,X_i](t)=X_i^{(2)}(t)$. The
processes $X_i^{(p_i)}=\{X_i^{(p_i)}(t),t\geq 0\}$,
$p_i=1,2,\cdots$, are again L\'{e}vy processes. They jump at the
same points as the original L\'{e}vy processes.

We have $E[X_i(t)]=E[X_i^{(1)}(t)]=tm_{i;1}<\infty$ and by Protter
(1990,p29), that
\begin{eqnarray}
E[X_i^{(p_i)}(t)]=E\left[\sum\limits_{0<s\leq t}(\triangle
X_i)^{p_i}\right]=t\int
x_i^{p_i}\nu(d\bm{x})=m_{i;p_i}t<\infty,\quad p_i\geq 2.\nonumber
\end{eqnarray}
Therefore, we can denote by
\begin{eqnarray}
Y_i^{(p_i)}(t)\stackrel{\rm{def}}{=}X_i^{(p_i)}(t)-E[X_i^{(p_i)}(t)]=X_i^{(p_i)}(t)-m_{i;p_i}t,\quad
p_i=1,2,3,\cdots\nonumber
\end{eqnarray}
the compensated power jump process of order $p_i$, and $Y_i^{(p_i)}$
is also a normal martingale.

We will express
$(X_1(t+t_0)-X_1(t_0))^{k_1}(X_2(t+t_0)-X_2(t_0))^{k_2}\cdots(X_n(t+t_0)-X_n(t_0))^{k_n}$,
$t_0$,$t\geq 0$, $k_i=1,2,\cdots$, $i=1,2,\cdots,n$, as a sum of
stochastic integrals with respect to the special processes
$Y_i^{(p_i)}(t)$, $i=1,\cdots,n$, $p_i=1,\cdots,k_i$.

Using It\^{o}'s formula(Protter, 1990, p.74, Theorem)we can write
for $k_i\geq 2$, $i=1,2,\cdots,n$,
\begin{eqnarray}
\begin{array}{rl}
&(X_1(t+t_0)-X_1(t_0))^{k_1}(X_2(t+t_0)-X_2(t_0))^{k_2}\cdots(X_n(t+t_0)-X_n(t_0))^{k_n}\\
=&\sum\limits_{i=1}^n\int_0^tk_i(X_i(s+t_0)-X_i(t_0))^{k_i-1}\prod\limits_{j\neq
i}(X_j(s+t_0)-X_j(t_0))dX_i(s)\\
&+\frac{1}{2}\sum\limits_{i=1}^n\int_0^t\sigma_{ii}^2k_i(k_i-1)(X_i(s+t_0)-X_i(t_0))^{k_i-2}\prod\limits_{j\neq
i}(X_j(s+t_0)-X_j(t_0))^{k_j}ds\\
&+\sum\limits_{1\leq i<j\leq
n}\int_0^t\sigma_{ij}^2k_ik_j(X_i(s+t_0)-X_i(t_0))^{k_i-1}(X_j(s+t_0)-X_j(t_0))^{k_j-1}\prod\limits_{\ell\neq
i,j}(X_{\ell}(s+t_0)-X_{\ell}(t_0))^{k_\ell}ds\\
&+\sum\limits_{0<s\leq
t}\left\{\sum\limits_{i=1}^n[(X_i(s+t_0)-X_i(t_0))^{k_i}-(X_i((s+t_0)-)-X_i(t_0))^{k_i}]
\prod\limits_{j\neq i}(X_j(s+t_0)-X_j(t_0))^{k_j}\right.\\
&\left.-\sum\limits_{i=1}^nk_i(X_i((s+t_0)-)-X_i(t_0))^{k_i-1}\prod\limits_{j\neq
i}(X_j(s+t_0)-X_j(t_0))^{k_j}\triangle X_i(s+t_0)\right\}
\end{array}\nonumber
\end{eqnarray}
\begin{eqnarray}
\begin{array}{rl}
=&\sum\limits_{i=1}^n\int_{t_0}^{t_0+t}k_i(X_i(u)-X_i(t_0))^{k_i-1}\prod\limits_{j\neq
i}(X_j(u)-X_j(t_0))dX_i^{(1)}(u)\\
&+\frac{1}{2}\sum\limits_{i=1}^n\sigma_{ii}^2k_i(k_i-1)\left[t(X_i(t+t_0)-X_i(t_0))^{k_i-2}\prod\limits_{j\neq
i}(X_j(t+t_0)-X_j(t_0))^{k_j}\right.\\
&\left.-\int_0^tsd\left((X_i(s+t_0)-X_i(t_0))^{k_i-2}\prod\limits_{j\neq
i}(X_j(s+t_0)-X_j(t_0))^{k_j}\right)\right]\\
&+\sum\limits_{1\leq i<j\leq
n}\sigma_{ij}^2k_ik_j\left[t(X_i(t+t_0)-X_i(t_0))^{k_i-1}(X_j(t+t_0)-X_j(t_0))^{k_j-1}\prod\limits_{\ell\neq
i,j}(X_{\ell}(t+t_0)-X_{\ell}(t_0))^{k_\ell}\right.\\
&-\left.\int_0^tsd\left
((X_i(s+t_0)-X_i(t_0))^{k_i-1}(X_j(s+t_0)-X_j(t_0))^{k_j-1}\prod\limits_{\ell\neq
i,j}(X_{\ell}(s+t_0)-X_{\ell}(t_0))^{k_\ell}ds\right)\right]\\
&+\sum\limits_{0<s\leq
t}\left\{\sum\limits_{i=1}^n[(X_i((s+t_0)-)+\triangle
X_i(s+t_0)-X_i(t_0))^{k_i}-(X_i((s+t_0)-)-X_i(t_0))^{k_i}]\right.\\
&\left.\times\prod\limits_{j\neq
i}(X_j(s+t_0)-X_j(t_0))^{k_j}\right.\\
&\left.-\sum\limits_{i=1}^nk_i(X_i((s+t_0)-)-X_i(t_0))^{k_i-1}\prod\limits_{j\neq
i}(X_j(s+t_0)-X_j(t_0))^{k_j}\triangle X_i(s+t_0)\right\}
\end{array}\nonumber
\end{eqnarray}

\begin{eqnarray}
\begin{array}{rl}
=&\sum\limits_{i=1}^n\int_{t_0}^{t_0+t}k_i(X_i(u)-X_i(t_0))^{k_i-1}\prod\limits_{j\neq
i}(X_j(u)-X_j(t_0))dX_i^{(1)}(u)\\
&+\frac{1}{2}\sum\limits_{i=1}^n\sigma_{ii}^2k_i(k_i-1)\left[t(X_i(t+t_0)-X_i(t_0))^{k_i-2}\prod\limits_{j\neq
i}(X_j(t+t_0)-X_j(t_0))^{k_j}\right.\\
&\left.-\int_0^tsd\left((X_i(s+t_0)-X_i(t_0))^{k_i-2}\prod\limits_{j\neq
i}(X_j(s+t_0)-X_j(t_0))^{k_j}\right)\right]\\
&+\sum\limits_{1\leq i<j\leq
n}\sigma_{ij}^2k_ik_j\left[t(X_i(t+t_0)-X_i(t_0))^{k_i-1}(X_j(t+t_0)-X_j(t_0))^{k_j-1}\prod\limits_{\ell\neq
i,j}(X_{\ell}(t+t_0)-X_{\ell}(t_0))^{k_\ell}\right.\\
&-\left.\int_0^tsd\left
((X_i(s+t_0)-X_i(t_0))^{k_i-1}(X_j(s+t_0)-X_j(t_0))^{k_j-1}\prod\limits_{\ell\neq
i,j}(X_{\ell}(s+t_0)-X_{\ell}(t_0))^{k_\ell}ds\right)\right]\\
&+\sum\limits_{0<s\leq
t}\sum\limits_{i=1}^n\sum\limits_{\ell=2}^{k_i}\left(
                                                 \begin{array}{c}
                                                   k_i \\
                                                   \ell \\
                                                 \end{array}
                                               \right)
(X_i((s+t_0)-)-X_i(t_0))^{k_i-\ell}\prod\limits_{j\neq
i}(X_j(s+t_0)-X_j(t_0))^{k_j}\left(\triangle X_i(s+t_0)\right)^\ell
\end{array}\nonumber
\end{eqnarray}

\begin{eqnarray}
\begin{array}{rl}
=&\sum\limits_{i=1}^n\int_{t_0}^{t_0+t}k_i(X_i(u)-X_i(t_0))^{k_i-1}\prod\limits_{j\neq
i}(X_j(u)-X_j(t_0))dX_i^{(1)}(u)\\
&+\frac{1}{2}\sum\limits_{i=1}^n\sigma_{ii}^2k_i(k_i-1)\left[t(X_i(t+t_0)-X_i(t_0))^{k_i-2}\prod\limits_{j\neq
i}(X_j(t+t_0)-X_j(t_0))^{k_j}\right.\\
&\left.-\int_0^tsd\left((X_i(s+t_0)-X_i(t_0))^{k_i-2}\prod\limits_{j\neq
i}(X_j(s+t_0)-X_j(t_0))^{k_j}\right)\right]\\
&+\sum\limits_{1\leq i<j\leq
n}\sigma_{ij}^2k_ik_j\left[t(X_i(t+t_0)-X_i(t_0))^{k_i-1}(X_j(t+t_0)-X_j(t_0))^{k_j-1}\prod\limits_{\ell\neq
i,j}(X_{\ell}(t+t_0)-X_{\ell}(t_0))^{k_\ell}\right.\\
&-\left.\int_0^tsd\left
((X_i(s+t_0)-X_i(t_0))^{k_i-1}(X_j(s+t_0)-X_j(t_0))^{k_j-1}\prod\limits_{\ell\neq
i,j}(X_{\ell}(s+t_0)-X_{\ell}(t_0))^{k_\ell}ds\right)\right]\\
&+\sum\limits_{t_0<u\leq
t+t_0}\sum\limits_{i=1}^n\sum\limits_{\ell=2}^{k_i}\left(
                                                 \begin{array}{c}
                                                   k_i \\
                                                   \ell \\
                                                 \end{array}
                                               \right)
\left(\triangle
X_i(u)\right)^\ell(X_i(u-)-X_i(t_0))^{k_i-\ell}\prod\limits_{j\neq
i}(X_j(u-)-X_j(t_0))^{k_j}\\
\end{array}\nonumber
\end{eqnarray}
\begin{eqnarray}
\begin{array}{rl}
=&\sum\limits_{i=1}^n\sum\limits_{\ell=2}^{k_i}\left(
                                                 \begin{array}{c}
                                                   k_i \\
                                                   \ell \\
                                                 \end{array}
                                               \right)
\int_{t_0}^{t+t_0}(X_i(u-)-X_i(t_0))^{k_i-\ell}\prod\limits_{j\neq
i}(X_j(u-)-X_j(t_0))^{k_j}dX_i^{(\ell)}(u)\\
&+\frac{1}{2}\sum\limits_{i=1}^n\sigma_{ii}^2k_i(k_i-1)\left[t(X_i(t+t_0)-X_i(t_0))^{k_i-2}\prod\limits_{j\neq
i}(X_j(t+t_0)-X_j(t_0))^{k_j}\right.\\
&\left.-\int_0^tsd\left((X_i(s+t_0)-X_i(t_0))^{k_i-2}\prod\limits_{j\neq
i}(X_j(s+t_0)-X_j(t_0))^{k_j}\right)\right]\\
&+\sum\limits_{1\leq i<j\leq
n}\sigma_{ij}^2k_ik_j\left[t(X_i(t+t_0)-X_i(t_0))^{k_i-1}(X_j(t+t_0)-X_j(t_0))^{k_j-1}\prod\limits_{\ell\neq
i,j}(X_{\ell}(t+t_0)-X_{\ell}(t_0))^{k_\ell}\right.\\
&-\left.\int_0^tsd\left
((X_i(s+t_0)-X_i(t_0))^{k_i-1}(X_j(s+t_0)-X_j(t_0))^{k_j-1}\prod\limits_{\ell\neq
i,j}(X_{\ell}(s+t_0)-X_{\ell}(t_0))^{k_\ell}ds\right)\right]
\end{array}\nonumber\\
\end{eqnarray}

\begin{lem}
The power of an increment of a L\'{e}vy process,
$(X_1(t+t_0)-X_1(t_0))^{k_1}(X_2(t+t_0)-X_2(t_0))^{k_2}\cdots(X_n(t+t_0)-X_n(t_0))^{k_n}$,
has a representation of the form
\begin{eqnarray}
\begin{array}{rl}
&(X_1(t+t_0)-X_1(t_0))^{k_1}(X_2(t+t_0)-X_2(t_0))^{k_2}\cdots(X_n(t+t_0)-X_n(t_0))^{k_n}\\
=&f^{(k)}(t,t_0)\\
&+\sum\limits_{m=1}^n\sum\limits_{q_{i_1}=1}^{k_1}\cdots\sum\limits_{q_{i_m}=1}^{k_m}\sum\limits_{\begin{array}{c}
                                                   \times_{j=1}^m(p_{i_j,1},\cdots,p_{i_j,q_{i_j}})\in \\
                                                   \times_{j=1}^m\{1,\cdots,k_j\}^{q_{i_j}}
                                                 \end{array}
}\int_{t_0}^{t+t_0}\int_{t_0}^{t_{i_1,1}-}\cdots\int_{t_0}^{t_{i_1,q_{i_1}}-}\cdots\int_{t_0}^{t_{i_m,1}-}\cdots\\
&\int_{t_0}^{t_{i_m,q_{i_m}}-}f_{(p_{i_1,1},\cdots,p_{i_1,q_{i_1}};\cdots;p_{i_m,1},\cdots,p_{i_m,q_{i_m}})}^{(\bm{k})}
(t,t_0;t_{i_1,1},\cdots,t_{i_1,q_{i_1}};\cdots;t_{i_m,1},\cdots,t_{i_m,q_{i_m}})\\
& dY_{i_m}^{(p_{i_m,q_{i_m}})}(t_{i_m,q_{i_m}})\cdots
dY_{i_m}^{(p_{i_m,2})}(t_{i_m,2})dY_{i_m}^{(p_{i_m,1})}(t_{i_m,1})\\
&\cdots dY_{i_1}^{(p_{i_1,q_{i_1}})}(t_{i_1,q_{i_1}})\cdots
dY_{i_1}^{(p_{i_1,2})}(t_{i_1,2})dY_{i_1}^{(p_{i_1,1})}(t_{i_1,1})
\end{array}\nonumber\\
\end{eqnarray}
where the
$f_{(p_{i_1,1},\cdots,p_{i_1,q_{i_1}};\cdots;p_{i_m,1},\cdots,p_{i_m,q_{i_m}})}^{(\bm{k})}
(t,t_0;t_{i_1,1},\cdots,t_{i_1,q_{i_1}};\cdots;t_{i_m,1},\cdots,t_{i_m,q_{i_m}})$
are deterministic functions in
$L^2(\mathbb{R}_{+}^{q_{i_1}+\cdots+q_{i_m}})$, and
$\bm{k}=(k_1,\cdots,k_n)$. In addition, the index $m$ controls the
number of $Y_i(t)$, $i=1,\cdots,n$, chosen from
$Y_1(t),\cdots,Y_n(t)$. After the $m$ is fixed, $(i_1,\cdots,i_m)$
indicates an arbitrary subset of the integer set $(1,2,\cdots,n)$.
After the $(i_1,\cdots,i_m)$ is fixed, for $i_j\in(i_1,\cdots,i_m)$,
the $(t_{i_j,1},\cdots,t_{i_j,q_{i_j}})$ with double index
$(i_j,\cdot)$ indicates the time-points chosen from the time-points
$(t_1,\cdots,t_{k_j})$ corresponding to
$(Y_{i_j}(t_1),\cdots,Y_{i_j}(t_{k_j}))$. The meaning of power index
$p_{i_j,\cdot}$ is similar.
\end{lem}

\noindent\textbf{Proof} Representation (3) follows from (2), where
we bring in the right compensations, i.e. we can write
\begin{eqnarray}
\begin{array}{rl}
&\sum\limits_{i=1}^n\sum\limits_{\ell=1}^{k_i}\left(
                                                \begin{array}{c}
                                                  k_i \\
                                                  \ell \\
                                                \end{array}
                                              \right)
\int_{t_0}^{t+t_0}(X_i(s-)-X_i(t_0))^{k_i-\ell}
\prod\limits_{j\neq i}(X_j(s-)-X_j(t_0))^{k_j}dX_i^{(\ell)}(s)\\
=&\sum\limits_{i=1}^n\sum\limits_{\ell=1}^{k_i}\left(
                                                \begin{array}{c}
                                                  k_i \\
                                                  \ell \\
                                                \end{array}
                                              \right)
\int_{t_0}^{t+t_0}(X_i(s-)-X_i(t_0))^{k_i-\ell} \prod\limits_{j\neq
i}(X_j(s-)-X_j(t_0))^{k_j}dY_i^{(\ell)}(s)\\
&+\sum\limits_{i=1}^n\sum\limits_{\ell=1}^{k_i}\left(
                                                \begin{array}{c}
                                                  k_i \\
                                                  \ell \\
                                                \end{array}
                                              \right)
m_{i\ell}\int_{t_0}^{t+t_0}(X_i(s-)-X_i(t_0))^{k_i-\ell}
\prod\limits_{j\neq i}(X_j(s-)-X_j(t_0))^{k_j}ds\\
=&\sum\limits_{i=1}^n\sum\limits_{\ell=1}^{k_i}\left(
                                                \begin{array}{c}
                                                  k_i \\
                                                  \ell \\
                                                \end{array}
                                              \right)
\int_{t_0}^{t+t_0}(X_i(s-)-X_i(t_0))^{k_i-\ell} \prod\limits_{j\neq
i}(X_j(s-)-X_j(t_0))^{k_j}dY_i^{(\ell)}(s)\\
&+\sum\limits_{i=1}^n\sum\limits_{\ell=1}^{k_i-1}\left(
                                                \begin{array}{c}
                                                  k_i \\
                                                  \ell \\
                                                \end{array}
                                              \right)
m_{i\ell}t(X_i(t+t_0)-X_i(t_0))^{k_i-\ell} \prod\limits_{j\neq
i}(X_j(t+t_0)-X_j(t_0))^{k_j}ds\\
&-\sum\limits_{i=1}^n\sum\limits_{\ell=1}^{k_i-1}\left(
                                                \begin{array}{c}
                                                  k_i \\
                                                  \ell \\
                                                \end{array}
                                              \right)
m_{i\ell}\int_{t_0}^{t+t_0}sd\left((X_i(s-)-X_i(t_0))^{k_i-\ell}
\prod\limits_{j\neq i}(X_j(s-)-X_j(t_0))^{k_j}\right)\\
&+\sum\limits_{i=1}^nm_{ik_i}t .
\end{array}\nonumber\\
\end{eqnarray}

Combining (2) and (4) gives
\begin{eqnarray}
\begin{array}{rl}
&\prod\limits_{i=1}^n(X_i(t+t_0)-X_i(t_0))^{k_i}\\
=&\frac{1}{2}\sum\limits_{i=1}^n\sigma_{ii}^2k_i(k_i-1)\left[t(X_i(t+t_0)-X_i(t_0))^{k_i-2}\prod\limits_{j\neq
i}(X_j(t+t_0)-X_j(t_0))^{k_j}\right.\\
&\left.-\int_0^tsd\left((X_i(s+t_0)-X_i(t_0))^{k_i-2}\prod\limits_{j\neq
i}(X_j(s+t_0)-X_j(t_0))^{k_j}\right)\right]\\
&+\sum\limits_{1\leq i<j\leq
n}\sigma_{ij}^2k_ik_j\left[t(X_i(t+t_0)-X_i(t_0))^{k_i-1}(X_j(t+t_0)-X_j(t_0))^{k_j-1}\prod\limits_{\ell\neq
i,j}(X_{\ell}(t+t_0)-X_{\ell}(t_0))^{k_\ell}\right.\\
&-\left.\int_0^tsd\left
((X_i(s+t_0)-X_i(t_0))^{k_i-1}(X_j(s+t_0)-X_j(t_0))^{k_j-1}\prod\limits_{\ell\neq
i,j}(X_{\ell}(s+t_0)-X_{\ell}(t_0))^{k_\ell}ds\right)\right]\\
&+\sum\limits_{i=1}^n\sum\limits_{\ell=1}^{k_i}\left(
                                                \begin{array}{c}
                                                  k_i \\
                                                  \ell \\
                                                \end{array}
                                              \right)
\int_{t_0}^{t+t_0}(X_i(s-)-X_i(t_0))^{k_i-\ell} \prod\limits_{j\neq
i}(X_j(s-)-X_j(t_0))^{k_j}dY_i^{(\ell)}(s)\\
&+\sum\limits_{i=1}^n\sum\limits_{\ell=1}^{k_i-1}\left(
                                                \begin{array}{c}
                                                  k_i \\
                                                  \ell \\
                                                \end{array}
                                              \right)
m_{i\ell}t(X_i(t+t_0)-X_i(t_0))^{k_i-\ell} \prod\limits_{j\neq
i}(X_j(t+t_0)-X_j(t_0))^{k_j}ds\\
&-\sum\limits_{i=1}^n\sum\limits_{\ell=1}^{k_i-1}\left(
                                                \begin{array}{c}
                                                  k_i \\
                                                  \ell \\
                                                \end{array}
                                              \right)
m_{i\ell}\int_{t_0}^{t+t_0}sd\left((X_i(s-)-X_i(t_0))^{k_i-\ell}
\prod\limits_{j\neq i}(X_j(s-)-X_j(t_0))^{k_j}\right)\\
&+\sum\limits_{i=1}^nm_{ik_i}t .
\end{array}\nonumber\\
\end{eqnarray}
The last equation is in terms of powers of increments of $X_i$ which
are strictly lower than $k_i$. So by induction representation (3)
can be proved. $\Box$

Notice that taking the expectation in (3) yields
\begin{eqnarray}
E\left[\prod\limits_{i=1}^n(X_i(t+t_0)-X_i(t_0))^{k_i}\right]
=f^{(\bm{k})}(t,t_0)=f^{(\bm{k})}(t),\quad t,t_0\geq 0,\nonumber
\end{eqnarray}
which is independent of $t_0$.

Moreover, it can easily be seen that
\begin{eqnarray}
f_{(p_{i_1,1},\cdots,p_{i_1,q_{i_1}};\cdots;p_{i_m,1},\cdots,p_{i_m,q_{i_m}})}^{(\bm{k})}
(t,t_0;t_{i_1,1},\cdots,t_{i_1,q_{i_1}};\cdots;t_{i_m,1},\cdots,t_{i_m,q_{i_m}})\nonumber
\end{eqnarray}
are just real multivariate polynomials of degree less than $k_i$ and
that we have
\begin{eqnarray}
f_{(p_{i_1,1},\cdots,p_{i_1,q_{i_1}};\cdots;p_{i_m,1},\cdots,p_{i_m,q_{i_m}})}^{(\bm{k})}
(t,t_0;t_{i_1,1},\cdots,t_{i_1,q_{i_1}};\cdots;t_{i_m,1},\cdots,t_{i_m,q_{i_m}})=0,\nonumber
\end{eqnarray}
whenever $p_{i_1,q_{i_1}}+\cdots+p_{i_j,q_{i_j}}>k_j$.

Because we can switch by a linear transformation from the
$Y_{i_m}^{(p_{i_m,2})}(t_{i_m,2})$ to
$H_{i_m}^{(p_{i_m,2})}(t_{i_m,2})$, it is clear that we also proved
the next representation.
\begin{lem}
The power of an increment of a L\'{e}vy process,
$(X_1(t+t_0)-X_1(t_0))^{k_1}(X_2(t+t_0)-X_2(t_0))^{k_2}\cdots(X_n(t+t_0)-X_n(t_0))^{k_n}$,
has a representation of the form
\begin{eqnarray}
\begin{array}{rl}
&(X_1(t+t_0)-X_1(t_0))^{k_1}(X_2(t+t_0)-X_2(t_0))^{k_2}\cdots(X_n(t+t_0)-X_n(t_0))^{k_n}\\
=&f^{(\bm{k})}(t,t_0)+\sum\limits_{d=1}^{|\bm{k}|}\sum\limits_{\bm{p}_1\in\mathbb{N}_d^n}
\int_{t_0}^{t+t_0}h_{(\bm{p}_1)}^{(\bm{k})}(t,t_0;t_1)dH^{\bm{p}_1}(t_1)\\
&+\sum\limits_{d=1}^{|\bm{k}|}\sum\limits_{\bm{p}_1+\bm{p}_2\in\mathbb{N}_d^n}
\int_{t_0}^{t+t_0}\int_{t_0}^{t_1}h_{(\bm{p}_1,\bm{p}_2)}^{(\bm{k})}(t,t_0;t_1,t_2)dH^{\bm{p}_2}(t_2)dH^{\bm{p}_1}(t_1)\\
&+\cdots\\
&+\sum\limits_{d=1}^{|\bm{k}|}\sum\limits_{\bm{p}_1+\cdots+\bm{p}_{|\bm{k}|}\in\mathbb{N}_d^n}
\int_{t_0}^{t+t_0}\int_{t_0}^{t_1}\cdots\int_{t_0}^{t_{|\bm{k}|-1}}
h_{(\bm{p}_1,\cdots,\bm{p}_{|\bm{k}|})}^{(\bm{k})}(t,t_0;t_1,t_2,\cdots,t_{|\bm{k}|})dH^{\bm{p}_N}(t_{|\bm{k}|})\cdots dH^{\bm{p}_2}(t_2)dH^{\bm{p}_1}(t_1)\\
=&\sum\limits_{m=1}^{|\bm{k}|}\sum\limits_{d=1}^{|\bm{k}|}\sum\limits_{\bm{p}_1+\cdots+\bm{p}_m\in\mathbb{N}_d^n}
\int_{t_0}^{t+t_0}\int_{t_0}^{t_1}\cdots\int_{t_0}^{t_{m-1}}
h_{(\bm{p}_1,\cdots,\bm{p}_m)}^{(\bm{k})}(t,t_0;t_1,t_2,\cdots,t_m)dH^{\bm{p}_m}(t_m)\cdots
dH^{\bm{p}_2}(t_2)dH^{\bm{p}_1}(t_1)
\end{array}\nonumber
\end{eqnarray}
where the index $m$ controls the number of integral variables in
each multiple integral,
$\bm{p}_i=(p_{i_1},\cdots,p_{i_n})\in\mathbb{N}_0^n$,
$\mathbb{N}_d^n=\{\bm{p}\in\mathbb{N}_0^n:|\bm{p}|=d\}$ and
$h_{(\bm{p}_1,\cdots,\bm{p}_m)}^{(\bm{k})}(t,t_0;t_1,t_2,\cdots,t_m)$
are deterministic functions in $L^2(\mathbb{R}_{+}^{m})$. Other
notations are the same as Lemma 1.
\end{lem}

\subsection{Representation of a square integrable random variable}
We first recall that $\{H_i^{(p_i)},i=1,2,\cdots,n;
p_i=1,2,\cdots\}$ is a set of pairwise strongly orthogonal
martingales, obtained by the orthogonalization procedure described
at the end of Section 2.

We denote by
\begin{eqnarray}
\begin{array}{l}
\mathscr{H}^{(\bm{p}_1,\cdots,\bm{p}_m)}=\left\{F\in L^2(\Omega):\right.\\
\int_{t_0}^{t+t_0}\int_{t_0}^{t_1}\cdots\int_{t_0}^{t_{m-1}}
h_{(\bm{p}_1,\cdots,\bm{p}_m)}^{(\bm{k})}(t,t_0;t_1,t_2,\cdots,t_m)dH^{\bm{p}_m}(t_m)\cdots
dH^{\bm{p}_2}(t_2)dH^{\bm{p}_1}(t_1)\\
\left.\bm{p}_j\in\mathbb{N}^n,\quad j=1,2,\cdots,m\right\}
\end{array}
\end{eqnarray}

We say that two multi-indexes
\begin{eqnarray}
(\bm{p}_1,\cdots,\bm{p}_m)\quad and \quad
(\tilde{\bm{p}}_1,\cdots,\tilde{\bm{p}}_{\tilde{m}}) \nonumber
\end{eqnarray}
are different if $m\neq \tilde{m}$ or when $m\equiv\tilde{m}$, if
there exists a subindex $1\leq l\leq m=\tilde{m}$, such that
$\bm{p}_l\neq \tilde{\bm{p}}_l$, and denote this by
\begin{eqnarray}
(\bm{p}_1,\cdots,\bm{p}_m)\neq
(\tilde{\bm{p}}_1,\cdots,\tilde{\bm{p}}_{\tilde{m}}) \nonumber
\end{eqnarray}

\begin{pro}
If
\begin{eqnarray}
(\bm{p}_1,\cdots,\bm{p}_m)\neq
(\tilde{\bm{p}}_1,\cdots,\tilde{\bm{p}}_{\tilde{m}}) \nonumber
\end{eqnarray}
then
\begin{eqnarray}
\mathscr{H}^{(\bm{p}_1,\cdots,\bm{p}_m)}\bot
\mathscr{H}^{(\tilde{\bm{p}}_1,\cdots,\tilde{\bm{p}}_{\tilde{m}})}.
\nonumber
\end{eqnarray}

\end{pro}
\noindent\textbf{Proof} Suppose we have two random variables $K\in
\mathscr{H}^{(\bm{p}_1,\cdots,\bm{p}_m)}$ and
$L\in\mathscr{H}^{(\tilde{\bm{p}}_1,\cdots,\tilde{\bm{p}}_{\tilde{m}})}$.
We need to prove that if
\begin{eqnarray}
(\bm{p}_1,\cdots,\bm{p}_m)\neq
(\tilde{\bm{p}}_1,\cdots,\tilde{\bm{p}}_{\tilde{m}}) \nonumber
\end{eqnarray}
then $K\bot L$.

For the case $m=\tilde{m}$, we use induction on $m$. Take first
$m=\tilde{m}=1$ and assume the following representations for $K$ and
$L$:
\begin{eqnarray}
K=\int_0^{\infty}f(t_1)dH^{\bm{p}_1}(t_1),\nonumber \quad
L=\int_0^{\infty}g(t_1)dH^{\tilde{\bm{p}}_1}(t_1)
\end{eqnarray}
where we must have $\bm{p}_1\neq \tilde{\bm{p}}_1$, By construction
$H^{\bm{p}_1}$ and $H^{\tilde{\bm{p}}_1}$ are strongly orthogonal
martingales. Using the fact that stochastic integrals with respect
to strongly orthogonal martingales are again strongly orthogonal
(Protter,1990,Lemma 2 and and Theorem 35, p.149) and thus also
weakly orthogonal, it immediately follows that $K\bot L$.

Suppose the theorem holds for all $1\leq m=\tilde{m}\leq
\textrm{n}-1$. We are going to prove the theorem for
$m=\tilde{m}=\textrm{n}$. Assume the following representations:
\begin{eqnarray}
\begin{array}{l}
K=\int_{t_0}^{t+t_0}\int_{t_0}^{t_1}\cdots\int_{t_0}^{t_{m-1}}
h_{(\bm{p}_1,\cdots,\bm{p}_m)}^{(\bm{k})}(t,t_0;t_1,t_2,\cdots,t_m)dH^{\bm{p}_m}(t_m)\cdots
dH^{\bm{p}_2}(t_2)dH^{\bm{p}_1}(t_1)\\
=\int_0^\infty\alpha(t_1)dH^{\bm{p}_1}(t_1)
\end{array}\nonumber
\end{eqnarray}

\begin{eqnarray}
\begin{array}{l}
L=\int_{t_0}^{t+t_0}\int_{t_0}^{t_1}\cdots\int_{t_0}^{t_{m-1}}
g_{(\tilde{\bm{p}}_1,\cdots,\tilde{\bm{p}}_m)}^{(\bm{k})}(t,t_0;t_1,t_2,\cdots,t_m)dH^{\tilde{\bm{p}}_m}(t_m)\cdots
dH^{\tilde{\bm{p}}_2}(t_2)dH^{\tilde{\bm{p}}_1}(t_1)\\
=\int_0^\infty\beta(t_1)dH^{\tilde{\bm{p}}_1}(t_1)\nonumber
\end{array}
\end{eqnarray}
There are two possibilities: (1) $\bm{p}_1=\tilde{\bm{p}}_1$ and (2)
$\bm{p}_1\neq\tilde{\bm{p}}_1$. In the former case we must have that
\begin{eqnarray}
(\bm{p}_2,\cdots,\bm{p}_m)\neq
(\tilde{\bm{p}}_2,\cdots,\tilde{\bm{p}}_{\tilde{m}}) \nonumber
\end{eqnarray}
 and thus by induction $\alpha(t_1)\bot \beta(t_1)$, so that
\begin{eqnarray}
E[KL]&=&E\left[\int_0^\infty\alpha_s\beta_sd<H^{\bm{p}_1},H^{\bm{p}_1}>_s\right]\nonumber\\
&=&\int_0^\infty
E(\alpha_s\beta_s)d<H^{\bm{p}_1},H^{\bm{p}_1}>_s=0.\nonumber
\end{eqnarray}
In the latter case we use again the fact that stochastic integrals
with respect to strongly orthogonal martingales are again strongly
orthogonal (Protter,1990,Lemma 2 and and Theorem 35, p.149) and thus
also weakly orthogonal. So it immediately follows that $K\bot L$.

For the case $m\neq \tilde{m}$, a similar argument can be used
together with the fact that all elements of every
$\mathscr{H}^{(\bm{p}_1,\cdots,\bm{p}_\ell)}$, $\ell\geq 1$, have
mean zero and thus are orthogonal w.r.t. the constants. $\Box$

\begin{pro}
Let
\begin{eqnarray}
\mathscr{P}=\left\{\prod\limits_{i=1}^n\prod\limits_{j=1}^m(X_i(t_j)-X_i(t_{j-1}))^{k_{i,j}}:
m\geq 0, 0\leq t_0\leq t_1<t_2<\cdots<t_m,
k_{1,1},\cdots,k_{n,m}\geq 1\right\}\nonumber
\end{eqnarray}
then we have that $\mathscr{P}$ is a total family in
$L^2(\Omega,\mathscr{F})$, i.e. the linear subspace spanned by
$\mathscr{P}$ is dense in $L^2(\Omega,\mathscr{F})$.
\end{pro}
\noindent\textbf{Proof} Let $Z\in L^2(\Omega,\mathscr{F})$ and
$Z\bot\mathscr{P}$. For any given $\varepsilon>0$, there exists a
finite set $\{0<s_1<\cdots<s_m\}$ and a square integrable random
variable $Z_\varepsilon\in
L^2(\Omega,\sigma(X(s_1),X(s_2),\cdots,X(s_m)))$ such that
\begin{eqnarray}
E\left[\|Z-Z_\varepsilon\|^2\right]<\varepsilon .\nonumber
\end{eqnarray}
So there exists a Borel function $f$ such that
\begin{eqnarray}
Z_\varepsilon=f_\varepsilon(X(s_1),X(s_2)-X(s_1),\cdots,X(s_m)-X(s_{m-1})).\nonumber
\end{eqnarray}
Because the polynomials are dense in
$L^2(\mathbb{R}^n,\mathbb{P}\circ X(t)^{-1})$ for each $t>0$, we can
approximate $Z_\varepsilon$ by polynomials. Furthermore because
$Z\bot\mathscr{P}$, we have $E[ZZ_\varepsilon]=0$. Then
\begin{eqnarray}
E\left[\|Z\|^2\right]=E[Z\cdot(Z-Z_\varepsilon)]\leq
\sqrt{E[\|Z\|^2]E\left[\|Z-Z_\varepsilon\|^2\right]}\leq
\sqrt{\varepsilon E[\|Z\|^2]},\nonumber
\end{eqnarray}
and Letting $\varepsilon\rightarrow 0$ yields $Z=\bm{0}$ a.s. Thus
$\mathscr{P}$ is a total family in $L^2(\Omega,\mathscr{F})$. $\Box$

We are now in a position to prove our main theorem.

\begin{thm}
(Chaotic representation property (CRP)). Every random variable $F$
in $L^2(\Omega,\mathscr{F})$ has a representation of the form
\begin{eqnarray}
\begin{array}{rl}
F=&\mathbb{E}(F)+\sum\limits_{m=1}^{\infty}\sum\limits_{d=1}^{\infty}\sum\limits_{\bm{p}_1+\cdots+\bm{p}_m\in\mathbb{N}_d^n}
\int_0^\infty\int_0^{t_1}\cdots\\
&\int_0^{t_{m-1}}
f_{(\bm{p}_1,\cdots,\bm{p}_m)}^{(\bm{k})}(t,t_0;t_1,t_2,\cdots,t_m)dH^{\bm{p}_m}(t_m)\cdots
dH^{\bm{p}_2}(t_2)dH^{\bm{p}_1}(t_1)
\end{array}\nonumber
\end{eqnarray}
where the
$f_{(\bm{p}_1,\cdots,\bm{p}_m)}^{(\bm{k})}(t,t_0;t_1,t_2,\cdots,t_m)$'s
are functions in $L^2(\mathbb{R}_{+}^{m})$.
\end{thm}

\noindent\textbf{Proof} Because $\mathscr{P}$ is a total family in
$L^2(\Omega,\mathscr{F})$, it is sufficient to prove that every
element of $\mathscr{P}$ has a representation of the desired form.
This follows from the fact that $\mathscr{P}$ is build up from terms
of the from
$\prod\limits_{i=1}^n\prod\limits_{j=1}^m(X_i(t_j)-X_i(t_{j-1}))^{k_{i,j}}$,
wherein every term has on its turn a representation of the form (6),
and we can nicely combine two terms in the desired representation.
Indeed, we have for all $k_i,l_i\geq 1$, $i=1,\cdots,n$, and $0\leq
t<s\leq u<v$, that the product of
$\prod\limits_{i=1}^n(X_i(s)-X_i(t))^{k_i}(X_i(v)-X_i(u))^{l_i}$ is
a sum of products of the form $AB$ where
\begin{eqnarray}
\begin{array}{rl}
A=&\int_t^s\int_t^{t_1-}\cdots\int_t^{t_{m-1}-}
h_{(\bm{p}_1,\cdots,\bm{p}_m)}^{(\bm{k})}(s,t;t_1,t_2,\cdots,t_m)dH^{\bm{p}_m}(t_m)\cdots
dH^{\bm{p}_2}(t_2)dH^{\bm{p}_1}(t_1)\nonumber
\end{array}
\end{eqnarray}

and
\begin{eqnarray}
\begin{array}{rl}
B=&\int_u^v\int_u^{u_1-}\cdots\int_u^{u_{\tilde{m}-1}-}
h_{(\tilde{\bm{p}}_1,\cdots,\tilde{\bm{p}}_{\tilde{m}})}^{(\bm{L})}(v,u;u_1,u_2,\cdots,u_{\tilde{m}})
dH^{\tilde{\bm{p}}_{\tilde{m}}}(t_{\tilde{m}})\cdots
dH^{\tilde{\bm{p}}_2}(t_2)dH^{\tilde{\bm{p}}_1}(t_1)\nonumber
\end{array}
\end{eqnarray}
where $m$ and $\bar{m}$ are two integers.

We can write
\begin{eqnarray}
\begin{array}{rl}
AB=&\int_u^v\int_u^{u_1-}\cdots\int_u^{u_{\tilde{m}-1}-}\int_t^s\int_t^{t_1-}\cdots\int_t^{t_{m-1}-}
h_{(\tilde{\bm{p}}_1,\cdots,\tilde{\bm{p}}_{\tilde{m}})}^{(\bm{L})}(v,u;u_1,u_2,\cdots,u_{\tilde{m}})\\
&h_{(\bm{p}_1,\cdots,\bm{p}_m)}^{(\bm{k})}(s,t;t_1,t_2,\cdots,t_m)dH^{\bm{p}_m}(t_m)\cdots
dH^{\bm{p}_2}(t_2)dH^{\bm{p}_1}(t_1)dH^{\tilde{\bm{p}}_{\tilde{m}}}(t_{\tilde{m}})\cdots
dH^{\tilde{\bm{p}}_2}(t_2)dH^{\tilde{\bm{p}}_1}(t_1)\\
=&\int_0^\infty\int_0^{u_1-}\cdots\int_0^{u_{\tilde{m}-1}-}\int_0^{u_{\tilde{m}}-}\int_0^{t_1-}\cdots\int_0^{t_{m-1}-}
1_{(u,v]}(u_1)1_{(u,u_1]}(u_2)\cdots 1_{(u,u_{\tilde{m}-1}]}(u_{\tilde{m}})\nonumber\\
&1_{(t,s]}(t_1)1_{(t,t_1]}(t_2)\cdots
1_{(t,t_{m-1}]}(t_m)h_{(\tilde{\bm{p}}_1,\cdots,\tilde{\bm{p}}_{\tilde{m}})}^{(\bm{L})}
(v,u;u_1,u_2,\cdots,u_{\tilde{m}})h_{(\bm{p}_1,\cdots,\bm{p}_m)}^{(\bm{k})}(s,t;t_1,t_2,\cdots,t_m)\\
&dH^{\bm{p}_m}(t_m)\cdots
dH^{\bm{p}_2}(t_2)dH^{\bm{p}_1}(t_1)dH^{\tilde{\bm{p}}_{\tilde{m}}}(t_{\tilde{m}})\cdots
dH^{\tilde{\bm{p}}_2}(t_2)dH^{\tilde{\bm{p}}_1}(t_1)
\end{array}\nonumber
\end{eqnarray}
and the desired representation follows. $\Box$

\begin{thm}
(Predictable representation property (PRP)). Every random variable
$F$ in $L^2(\Omega,\mathscr{F})$ has a representation of the form
\begin{eqnarray}
\begin{array}{rl}
F=&\mathbb{E}(F)+\sum\limits_{d=1}^{\infty}\sum\limits_{\bm{p}\in\mathbb{N}_d^n}
\int_0^\infty \Phi^{\bm{p}}(s)dH^{\bm{p}}(t_m)(s)
\end{array}\nonumber
\end{eqnarray}
where $\Phi^{\bm{p}}(s)$ is predictable.
\end{thm}

\noindent\textbf{Proof} From the above theorem, we know that $F$ has
a representation of the form
\begin{eqnarray}
\begin{array}{rl}
&F-\mathbb{E}(F)\\
=&\sum\limits_{m=1}^{\infty}\sum\limits_{d=1}^{\infty}\sum\limits_{\bm{p}_1+\cdots+\bm{p}_m\in\mathbb{N}_d^n}
\int_0^\infty\int_0^{t_1-}\cdots\int_0^{t_{m-1}-}
f_{(\bm{p}_1,\cdots,\bm{p}_m)}^{(\bm{k})}(t_1,t_2,\cdots,t_m)dH^{\bm{p}_m}(t_m)\cdots
dH^{\bm{p}_2}(t_2)dH^{\bm{p}_1}(t_1)\\
=&\sum\limits_{d=1}^{\infty}\sum\limits_{\bm{p}_1\in\mathbb{N}_d^n}
\int_0^\infty f_{(\bm{p}_1)}^{(\bm{k})}(t_1)dH^{\bm{p}_1}(t_1)
+\sum\limits_{d=1}^{\infty}\sum\limits_{\bm{p}_1\in\mathbb{N}_d^n}
\int_0^\infty\left[\sum\limits_{k=0}^\infty\sum\limits_{\bm{p}_2+\cdots+\bm{p}_m\in\mathbb{N}_{k}^n}
\int_0^{t_1-}\cdots\right.\\
&\left. \int_0^{t_{m-1}}
f_{(\bm{p}_1,\cdots,\bm{p}_m)}^{(\bm{k})}(t_1,t_2,\cdots,t_m)dH^{\bm{p}_m}(t_m)\cdots
dH^{\bm{p}_2}(t_2)\right]dH^{\bm{p}_1}(t_1)\\
=&\sum\limits_{d=1}^{\infty}\sum\limits_{\bm{p}_1\in\mathbb{N}_d^n}
\int_0^\infty \left[f_{(\bm{p}_1)}^{(\bm{k})}(t_1)\right.\\
&\left.+\sum\limits_{k=0}^\infty\sum\limits_{\bm{p}_2+\cdots+\bm{p}_m\in\mathbb{N}_{k}^n}
\int_0^{t_1-}\int_0^{t_{m-1}}
f_{(\bm{p}_1,\cdots,\bm{p}_m)}^{(\bm{k})}(t_1,t_2,\cdots,t_m)dH^{\bm{p}_m}(t_m)\cdots
dH^{\bm{p}_2}(t_2)\right]dH^{\bm{p}_1}(t_1)\\
=&\sum\limits_{d=1}^{\infty}\sum\limits_{\bm{p}\in\mathbb{N}_d^n}
\int_0^\infty \Phi^{\bm{p}}(s)dH^{\bm{p}}(t_m)(s)
\end{array}\nonumber
\end{eqnarray}
which is exactly of the form we want. $\Box$

\begin{rmk}
Because we can identify every martingale $M\in\mathscr{U}^2$ with
its terminal value $M_{\infty}\in L^2(\Omega,\mathscr{F})$ and
because $M_t=\mathbb{E}[M_\infty|\mathscr{F}_t]$, we have the
predictable representation
\begin{eqnarray}
M_t&=&\sum\limits_{d=1}^{\infty}\sum\limits_{\bm{p}\in\mathbb{N}_d^n}
\int_0^t\Phi^{\bm{p}}(s)dH^{\bm{p}}(t_m)(s) \nonumber
\end{eqnarray}
which is a sum of strongly orthogonal martingales.
\end{rmk}

Another consequence of the chaotic representation property, is the
following theorem:
\begin{thm}
We have the following space decomposition:
\begin{eqnarray}
L^2(\Omega,\mathscr{F})=\mathbb{R}\oplus\left(\bigoplus_{d=1}^\infty\bigoplus_{\bm{p}\in\mathbb{N}_d^n}
\mathscr{H}^{\bm{p}}\right).\nonumber
\end{eqnarray}
\end{thm}

\begin{rmk}
The L\'{e}vy-Khintchine formula has a simpler expression when the
sample paths of the related L\'{e}vy process have bounded variation
on every compact time interval a.s. It is well known(Bertoin,
1996,p.15), that a L\'{e}vy process has bounded variation if and
only if $\Sigma=\bm{0}$, and
$\int(1\wedge\|\bm{x}\|)\nu(d\bm{x})<\infty$. In that case the
characteristic exponent can be re-expressed as
\begin{eqnarray}
\psi(\bm{\theta})=i\bm{d}\cdot
\bm{\theta}+\int_{\mathbb{R}^n}\left(exp(\textrm{i}\bm{\theta}\cdot
\bm{x})-1\right)\nu(d\bm{x}).\nonumber
\end{eqnarray}
Furthermore, we can write
\begin{eqnarray}
X_i(t)=dt+\sum\limits_{0<s\leq t}\triangle X_i(s),\quad t\geq
0,\quad i=1,2,\cdots,n.
\end{eqnarray}
and the calculations simplify somewhat because $\Sigma=\bm{0}$ and
for $k\geq 1$,
\begin{eqnarray}
\begin{array}{rl}
&\sum\limits_{i=1}^n\int_0^tk_i(X_i(s+t_0)-X_i(t_0))^{k_i-1}\prod\limits_{j\neq
i}(X_j(s+t_0)-X_j(t_0))dX_i(s)\nonumber\\
=&\sum\limits_{0<s\leq
t}\left\{\sum\limits_{i=1}^nk_i(X_i((s+t_0)-)-X_i(t_0))^{k_i-1}\prod\limits_{j\neq
i}(X_j(s+t_0)-X_j(t_0))^{k_j}\triangle X_i(s+t_0)\right\}\nonumber
\end{array}
\end{eqnarray}
\end{rmk}

\section{Examples}
\label{sec4}

Multidimensional models with jumps are more difficult to construct
than one-dimensional ones. A simple method to introduce jumps into a
multidimensional model is to take a multivariate Brownian motion and
time change it with a univariate subordinator (refer to Cont and
Tankov (2004)). The multidimensional versions of the models include
variance gamma, normal inverse Gaussian and generalized hyperbolic
processes. The principal advantage of this method is its simplicity
and analytic tractability; in particular, processes of this type are
easy to simulate. Another method to introduce jumps into a
multidimensional model is so-called method of L\'{e}vy copulas
proposed by Kallsen and Tankov (2006). The principle advantage in
this way lies in that the dependence among components of the
multidimensional L\'{e}vy processes can be completely characterized
with a L\'{e}vy copula. This allows us to give a systematic method
to construct multidimensional L\'{e}vy processes with specified
dependence.

In the following first and third examples, we define a multivariate
gamma process and a multivariate Meixner process by using L\'{e}vy
copulas, and furthermore discuss their orthogonalization procedures.
All the concepts and notations are adopted from the Kallsen and
Tankov (2006). In particular, for $n\geq 2$, the L\'{e}vy copula
$F(u_1,\cdots,u_n): \bar{\mathbb{R}}^n\rightarrow\bar{\mathbb{R}}$
is taken as
\begin{eqnarray}
F(u_1,\cdots,u_n)=2^{2-n}\left(\sum\limits_{j=1}^n|u_j|^{-\theta}\right)^{-1/\theta}(\eta
I_{\{u_1\cdots u_n\geq 0\}}-(1-\eta)I_{\{u_1\cdots u_n< 0\}}) .
\end{eqnarray}
It defines a two parameter family of L\'{e}vy copulas which
resembles the Clayton family of ordinary copulas. It is in fact a
L\'{e}vy copula homogeneous of order 1, for any $\theta>0$ and any
$\eta\in [0,1]$.

In addition, we know that if the tail integrals $U_i(x_i)$,
$i=1,\cdots,n$, are absolutely continuous, we can compute the
L\'{e}vy density of the L\'{e}vy copula process by differentiation
as follows:
\begin{eqnarray}
\nu(dx_1,\cdots,dx_n)=\partial_1\cdots\partial_nF|_{\xi_1=U_1(x_1),\cdots,\xi_n=U_n(x_n)}\nu_1(dx_1)\cdots\nu_1(dx_n)
\end{eqnarray}
where $\nu_1(dx_1),\cdots,\nu_n(x_n)$ are marginal L\'{e}vy
densities.

\subsection{The multivariate gamma process}
In the literature, the multivariate gamma distributions on
$\mathbb{R}^n$ have several non-equivalent definitions(refer to
Johnson and Balakrishnan (1997)). Here we consider only a
multivariate gamma process by using copula. The multivariate
\textsl{Gamma process} $\bm{G}(t)=(G_1(t),G_2(t),\cdots,G_n(t))^T$
is a multivariate L\'{e}vy process with the marginal distribution
density functions of $G_i(t)$, $i=1,2,\cdots,n$ given by
\begin{eqnarray}
f_{G_i(t)}(x_i)&=&\frac{1}{\Gamma(\gamma_it)}
\lambda_i^{\gamma_it}x_i^{\gamma_it-1}
exp\left\{-\lambda_ix_i\right\},\nonumber\\
&&x_i>0,\quad \lambda_i,\gamma_i>0\quad i=1,2,\cdots,n.\nonumber
\end{eqnarray}
The corresponding marginal characteristic functions are given by
\begin{eqnarray}
\mathbb{E}(e^{\theta_iG_i(t)})=\left(1-\frac{\textrm{i}\theta_i}{\lambda_i}\right)^{-\gamma_it},\qquad
i=1,2,\cdots, n.\nonumber
\end{eqnarray}
The corresponding marginal L\'{e}vy measures are given by
\begin{eqnarray}
\nu_i(dx_i)=&\frac{\gamma_i}{x_i}exp\{-\lambda_ix_i\}
I_{(0,\infty)}(x_i),\quad i=1,2,\cdots,n\nonumber
\end{eqnarray}

The corresponding L\'{e}vy measure is given by
\begin{eqnarray}
\begin{array}{rl}
\nu_P(d\bm{x})=&\partial_1\cdots\partial_nF|_{\xi_1=U_1(x_1),\cdots,\xi_n=U_n(x_n)}\frac{\prod\limits_{i=1}^n\gamma_i}
{\prod\limits_{i=1}^nx_i}exp\{-\sum\limits_{i=1}^n\lambda_ix_i\}
I_{(0,\infty)^n}(x_1,\cdots,x_n),\nonumber
\end{array}
\end{eqnarray}
where $d\bm{x}=dx_1\cdots dx_n $. The n-dimensional Gamma processes
are used i.a. in insurance mathematics(Dickson and Waters, 1993,
1996; Dufresne and Gerber, 1993;Dufresne et al., 1991).

We denote by
\begin{eqnarray}
G_i^{(p_i)}(t)=\sum\limits_{0<s\leq t}(\triangle G_i(s))^{p_i},\quad
p_i\geq 1,\quad i=1,2,\cdots,n \nonumber
\end{eqnarray}
the power jump processes of $G_i(t)$. In addition, set
$\bm{G}^{(p_1,\cdots,p_n)}(t)=(G_1^{(p_1)},\cdots,G_n^{(p_n)})$,
where $(p_1,\cdots,p_n)\in\mathbb{N}^n$. Using the exponential
formula (Bertoin,1996), and the change of the variable
$\bm{z}=(x_1^{p_1},\cdots, x_n^{p_n})$, we obtain for
$p_1+\cdots+p_n\geq 1$
\begin{eqnarray}
\begin{array}{l}
E\left[exp\left(\textrm{i}\bm{\theta}^T\bm{G}^{(p_1,\cdots,p_n)}(t)\right)\right]\\
=exp\left(t\int_{\mathbb{R}_{+}^n}
\left(exp(\textrm{i}\sum\limits_{i=1}^n\theta_ix_i^{p_i})-1\right)\partial_1\cdots\partial_nF|_{\xi_1=U_1(x_1),\cdots,\xi_n=U_n(x_n)}
\frac{\prod\limits_{i=1}^n\gamma_i}{\prod\limits_{i=1}^nx_i}
exp\{-\sum\limits_{i=1}^n\lambda_ix_i\}d\bm{x}\right)\\
=exp\left(t\int_{\mathbb{R}_{+}^n}
(exp(\textrm{i}\bm{\theta}^T\bm{z})-1)\partial_1\cdots\partial_nF|_{\xi_1=U_1(z_1^{1/p_1}),\cdots,\xi_n=U_d(z_n^{1/p_n})}
\frac{\prod\limits_{i=1}^n\gamma_i}{\prod\limits_{i=1}^np_iz_i}
exp\{-\sum\limits_{i=1}^n\lambda_iz_i^{1/p_i}\}d\bm{z}\right)
\end{array}\nonumber
\end{eqnarray}
which means that the L\'{e}vy measure of $\bm{G}^{(p_1,\cdots,p_n)}$
is
\begin{eqnarray}
\partial_1\cdots\partial_nF|_{\xi_1=U_1(z_1^{1/p_1}),\cdots,\xi_n=U_d(z_n^{1/p_n})}\frac{\prod\limits_{i=1}^n\gamma_i}{\prod\limits_{i=1}^np_iz_i}
exp\{-\sum\limits_{i=1}^n\lambda_iz_i^{1/p_i}\}d\verb"z" \nonumber
\end{eqnarray}

Introduce power jump processes of the form
\begin{eqnarray}
G(t)^{(p_1,\cdots,p_n)}\stackrel{\rm{def}}{=}\sum\limits_{0<s\leq
t}(\triangle G_1(s))^{p_1}\cdots (\triangle G_n(s))^{p_1}\nonumber
\end{eqnarray}
and then define the Teugels martingale monomial
\begin{eqnarray}
\hat{G}(t)^{(p_1,\cdots,p_n)}\stackrel{\rm{def}}{=}G(t)^{(p_1,\cdots,p_n)}-\mathbb{E}[G(t)^{(p_1,\cdots,p_n)}]
=G(t)^{(p_1,\cdots,p_n)}-m_{\bm{p}}t .\nonumber
\end{eqnarray}

Because
\begin{eqnarray}
\begin{array}{l}
E\left[\sum\limits_{0<s\leq t}(\triangle
G_1(s))^{p_1}\cdots(\triangle
G_n(s))^{p_n}\right]\\
=t\int_{\mathbb{R}_{+}^n}x_1^{p_1}\cdots
x_n^{p_n}\partial_1\cdots\partial_nF|_{\xi_1=U_1(x_1),\cdots,\xi_n=U_n(x_n)}\frac{\prod\limits_{i=1}^n\gamma_i}{\prod\limits_{i=1}^nx_i}
exp\{-\sum\limits_{i=1}^n\lambda_ix_i\}d\verb"x",\quad |\bm{p}|\geq
1,
\end{array}\nonumber
\end{eqnarray}

Next, we orthogonalize the set $\hat{G}^{\bm{p}}$ of martingales. So
we are looking for a set of martingales
\begin{eqnarray}
H^{\bm{p}}&=&\hat{G}^{\bm{p}}+\sum\limits_{\bm{q}\prec\bm{p},|\bm{q}|=|\bm{p}|}c_{\bm{q}}
\hat{G}^{\bm{q}}+\sum\limits_{k=1}^{|\bm{p}|-1}\sum\limits_{|\bm{q}|=k}c_{\bm{q}}\hat{G}^{\bm{q}},
\end{eqnarray}
such that $H^{\bm{p}}$ is strongly orthogonal to
$H^{\tilde{\bm{p}}}$, for $\bm{p}\neq \tilde{\bm{p}}$.

The first space $S_1$ in the gamma case is defined as follows
\begin{eqnarray}
S_1&=&\left\{\sum\limits_{k=1}^d\sum\limits_{|\bm{p}|=k}c_{k}(p_1,\cdots,p_n)x_1^{p_1}\cdots
x_n^{p_n}+c_0+\sum\limits_{(i_1,\cdots,i_n)\in\{0,-1\}^n,|\bm{i}|\geq
-(n-1)}c_{-1}(i_1,\cdots,i_n)x_1^{i_1}\cdots
x_n^{i_n}; \right.\nonumber\\
&&\left. d\in\{1,2,\cdots\},
c_{j}(p_1,\cdots,p_n)\in\mathbb{R},j=-1,0,\cdots,d;
x_i>0,i=1,\cdots,n;\bm{i}=(i_1,\cdots,i_n)\right\}\nonumber
\end{eqnarray}
which is endowed with a scalar product $<\cdot,\cdot>_1$, given by
\begin{eqnarray}
\begin{array}{rl}
&<P(\verb"x"),Q(\verb"x")>_1\nonumber\\
=&\int_{0}^{+\infty}\cdots\int_{0}^{+\infty}
P(\verb"x")Q(\verb"x")\partial_1\cdots\partial_nF|_{\xi_1=U_1(x_1),\cdots,\xi_n=U_n(x_n)}\prod\limits_{i=1}^n
(\gamma_ix_i)exp\{-\sum\limits_{i=1}^n\lambda_ix_i\}d\verb"x"
.\nonumber
\end{array}
\end{eqnarray}
Note that
\begin{eqnarray}
\begin{array}{l}
<x_1^{p_1-1}\cdots x_n^{p_n-1},x_1^{q_1-1}\cdots
x_n^{q_n-1}>_1\\
=\int_{\mathbb{R}_{+}^n}x_1^{p_1+q_1-1}\cdots
x_n^{p_n+q_n-1}\partial_1\cdots\partial_nF|_{\xi_1=U_1(x_1),\cdots,\xi_n=U_n(x_n)}\left(\prod\limits_{i=1}^n\gamma_i\right)
exp\{-\sum\limits_{i=1}^n\lambda_ix_i\}d\verb"x"\\
\quad |\bm{p}|,|\bm{q}|\geq 1. \nonumber
\end{array}
\end{eqnarray}

Thus we can construct the other space $S_2$ which is the space of
all linear transformations of the Teugels martingale monomials of
the multi-dimensional Gamma process, i.e.
\begin{eqnarray}
S_2&=&\left\{\sum\limits_{p_1+\cdots+p_n=d}a_{d}(p_1,\cdots,p_n)\hat{G}(t)^{(p_1,\cdots,p_n)}
+\sum\limits_{p_1+\cdots+p_n=d-1}a_{d-1}(p_1,\cdots,p_n)\hat{G}(t)^{(p_1,\cdots,p_n)}\right.\nonumber\\
&&\left.+\cdots+\sum\limits_{p_1+\cdots+p_n=1}a_{1}(p_1,\cdots,p_n)\hat{G}(t)^{(p_1,\cdots,p_n)},\quad
d\geq 1. \right\}.\nonumber
\end{eqnarray}
endowed with the scalar product $<\cdot,\cdot>_2$, given by
\begin{eqnarray}
\begin{array}{l}
<\hat{G}^{(p_1,\cdots,p_n)},\hat{G}^{(q_1,\cdots,q_n)}>_2\\
=E\left[\left[\hat{G}^{(p_1,\cdots,p_n)},\hat{G}^{(q_1,\cdots,q_n)}\right](1)\right]\\
=E\left[\hat{G}(1)^{(p_1+q_1,\cdots,p_n+q_n)}\right]\\
=\int_{\mathbb{R}_{+}^n}x_1^{p_1+q_1-1}\cdots
x_n^{p_n+q_n-1}\partial_1\cdots\partial_nF|_{\xi_1=U_1(x_1),\cdots,\xi_n=U_n(x_n)}\left(\prod\limits_{i=1}^n\gamma_i\right)
exp\{-\sum\limits_{i=1}^n\lambda_ix_i\}d\bm{x}
\end{array} \nonumber
\end{eqnarray}
So one clearly sees that $x_1^{p_1-1}x_2^{p_2-1}\cdots
x_n^{p_n-1}\leftrightarrow \hat{G}^{(p_1,\cdots,p_n)}$ is an
isometry between $S_{1}$ and $S_{2}$. An orthogonalization of
$\{x_1^{-1}x_2^{-1}\cdots x_{n-1}^{-1},x_1^{-1}x_3^{-1}\cdots
x_n^{-1},\cdots,x_n^{-1},1,x_1,\cdots,x_n,x_1^2,x_1x_2,\cdots,x_n^2,\cdots\}$
in $S_{1}$ can give the multivariate polynomials, so by isometry we
also can find an orthogonalization of
$\{\hat{G}^{(1,0,\cdots,0)},\cdots,\hat{G}^{(0,\cdots,0,1)}$,
$\hat{G}^{(2,0,\cdots,0)},\hat{G}^{(1,1,0,\cdots,0)},
\cdots,\hat{G}^{(0,\cdots,0,2)},\cdots\}$.

\subsection{The negative multinomial processes }

The next process of bounded variation we look at is the negative
multinomial processes, sometimes also called Pascal processes. Here
the conception of negative multinomial processes can be found in
Johnson et al.(1997), and P. Bernardoff (2003).

We define a negative multinomial distribution on $\mathbb{N}_0^n$.
Its distribution is
$\sum\limits_{\bm{k}\in\mathbb{N}_0^n}prob_{\bm{k}}\delta_{\bm{k}}$,
where $prob_{\bm{k}}\delta_{\bm{k}}$ denotes the probability measure
concentrated at $\bm{k}=\{k_1,k_2,\cdots,k_n\}$,
\begin{eqnarray}
prob_{\bm{k}}&\stackrel{\rm{def}}{=}&\mathbb{P}(k_1,\cdots,k_n)=\frac{\Gamma(t+\sum\limits_{i=1}^nk_i)}{k_1!k_2!\cdots
k_n!\Gamma(t)}\lambda^t\prod\limits_{i=1}^n(\mu\lambda_i)^{k_i},\nonumber\\
&&k_i=0,1,2,\cdots,\quad i=1,2,\cdots .\nonumber
\end{eqnarray}
where $0<\lambda<1$, $0<\mu\lambda_i<1$ for $i=1,2,\cdots,n$ and
$\lambda+\mu(\lambda_1+\cdots+\lambda_n)=1$.

This type of n-dimensional l\'{e}vy processes $P=\{P(t),t\geq 0\}$
where $P(t)=(P_1(t),P_2(t),\cdots,P_n(t))$, has a characteristic
function given by
\begin{eqnarray}
\mathbb{E}\left[exp(i\bm{\theta}'\cdot
P(t))\right]&=&\sum\limits_{\bm{k}\in \mathbb{N}^n}prob_{\bm{k}}
e^{ik_1\theta_1}\cdots e^{ik_n\theta_n}\nonumber\\
&=&\left(\frac{\lambda}{1-\mu(\lambda_1e^{i\theta_1}+\cdots+\lambda_ne^{i\theta_n})}\right)^t\nonumber
\end{eqnarray}

The corresponding L\'{e}vy measure $v(k_1,\cdots,k_n)$ is given by
\begin{eqnarray}
v(k_1,\cdots,k_n)=\frac{(|\bm{k}|-1)!}{k_1!\cdots
k_n!}\prod\limits_{i=1}^n\left(\mu\lambda_i\right)^{k_i},\qquad
|\bm{k}|\stackrel{\rm{def}}{=}k_1+\cdots+k_n.\nonumber
\end{eqnarray}
Let us denote with
\begin{eqnarray}
G(t)^{(p_1,\cdots,p_n)}\stackrel{\rm{def}}{=}\sum\limits_{0<s\leq
t}(\triangle P_1(s))^{p_1}\cdots (\triangle P_n(s))^{p_1},\qquad
|\bm{p}|\geq 1.\nonumber
\end{eqnarray}
the power jump processes of $P$ and with
$Q^{(p_1,\cdots,p_n)}=\{Q^{(p_1,\cdots,p_n)}(t),t\geq 0\}$ the
corresponding processes of Teugels martingale monomials.

We look for the orthogonalization of the set
$\{Q^{(p_1,\cdots,p_n)}, |\bm{p}|\geq 1\}$ of martingales. The space
$S_1$ is now defined as follows
\begin{eqnarray}
S_1&=&\left\{\sum\limits_{k=1}^d\sum\limits_{|\bm{p}|=k}c_{k}(p_1,\cdots,p_n)k_1^{p_1}\cdots
k_n^{p_n}+c_0+\sum\limits_{(i_1,\cdots,i_n)\in\{0,-1\}^n,|\bm{i}|\geq
-(n-1)}c_{-1}(i_1,\cdots,i_n)k_1^{i_1}\cdots
k_n^{i_n}; \right.\nonumber\\
&&\left. d\in\{1,2,\cdots\},
c_{j}(p_1,\cdots,p_n)\in\mathbb{R},j=-1,0,\cdots,d;
k_i\in\mathbb{N},i=1,\cdots,n;\bm{i}=(i_1,\cdots,i_n)\right\}\nonumber
\end{eqnarray}
endowed with a scalar product $<\cdot,\cdot>_1$, given by
\begin{eqnarray}
<P(\bm{k}),R(\bm{k})>_1=\sum\limits_{\bm{k}\in\mathbb{N}_0^n}
P(\bm{k})R(\bm{k})k_1\cdots k_n\frac{(|\bm{k}|-1)!}{k_1!\cdots
k_n!}\prod\limits_{i=1}^n\left(\mu\lambda_i\right)^{k_i}. \nonumber
\end{eqnarray}
Note that
\begin{eqnarray}
\begin{array}{rl}
<k_1^{p_1-1}\cdots k_n^{p_n-1},k_1^{q_1-1}\cdots
k_n^{q_n-1}>_1&=\sum\limits_{\bm{k}\in\mathbb{N}_0^n}k_1^{p_1+q_1-1}\cdots
k_n^{p_1+q_2-1}\frac{(|\bm{k}|-1)!}{k_1!\cdots
k_n!}\prod\limits_{i=1}^n\left(\mu\lambda_i\right)^{k_i}\\
&\quad p_i,q_j\geq 1,\quad i,j=1,2,\cdots,n \nonumber
\end{array}
\end{eqnarray}
The other space $S_2$ is the space of all linear transformations of
the Teugels martingales of the negative multinomial processes, i.e.
\begin{eqnarray}
S_2&=&\left\{\sum\limits_{p_1+\cdots+p_n=d}a_{d}(p_1,\cdots,p_n)Q(t)^{(p_1,\cdots,p_n)}
+\sum\limits_{p_1+\cdots+p_n=d-1}a_{d-1}(p_1,\cdots,p_n)Q(t)^{(p_1,\cdots,p_n)}\right.\nonumber\\
&&\left.\qquad+\cdots+\sum\limits_{p_1+\cdots+p_n=1}a_{1}(p_1,\cdots,p_n)Q(t)^{(p_1,\cdots,p_n)},\right.\nonumber\\
&&\quad d\in\mathbb{N}, \quad (p_1,\cdots,p_n)\in\mathbb{N}^n,\quad
a_{j}(p_1,\cdots,p_n)\in\mathbb{R},\quad j=1,2,\cdots,d \}.\nonumber
\end{eqnarray}
and is endowed with the scalar product $<\cdot,\cdot>_2$, given by
\begin{eqnarray}
<Q^{(p_1,\cdots,p_n)},Q^{(q_1,\cdots,q_n)}>_2&=&\mathbb{E}\left[\left[Q^{(p_1,\cdots,p_n)},Q^{(q_1,\cdots,q_n)}\right]_1\right]\nonumber\\
&=&\mathbb{E}\left[Q^{(p_1+q_1,\cdots,p_n+q_n)}(1)\right]\nonumber\\
&=&\sum\limits_{\bm{k}\in\mathbb{N}^n}k_1^{p_1+q_1}\cdots
k_n^{p_n+q_n}\frac{(|\bm{k}|-1)!}{k_1!\cdots
k_n!}\prod\limits_{i=1}^n\left(\mu\lambda_i\right)^{k_i} .\nonumber
\end{eqnarray}
So one clearly sees that $k_1^{p_1-1}k_2^{p_2-1}\cdots
k_n^{p_n-1}\leftrightarrow Q^{(p_1,\cdots,p_n)}$ is an isometry
between $S_{1}$ and $S_{2}$. An orthogonalization of
$\{k_1^{-1}k_2^{-1}\cdots k_{n-1}^{-1},k_1^{-1}k_3^{-1}\cdots
k_n^{-1},\cdots,k_n^{-1},1,k_1,\cdots,k_n,k_1^2,k_1k_2,\cdots,k_n^2,\cdots\}$
in $S_{2}$ gives the multivariate Meixner
polynomials(Griffiths(1975), Griffiths and Span\`{o}(2008), and
Koekoek and Swarttouw(1998)), so by isometry we also find an
orthogonalization of
$$\{Q^{(1,0,\cdots,0)},\cdots,Q^{(0,\cdots,0,1)},Q^{(2,0,\cdots,0)},Q^{(1,1,0,\cdots,0)},\cdots,G^{(0,\cdots,0,2)},\cdots\}$$.

\subsection{The Multivariate Meixner process}
A multivariate Meixner process
$M(t)=(M_1(t),M_2(t),\cdots,M_n(t))^T$,$t\geq 0$ is a bounded
variation L\'{e}vy process based on the infinitely divisible
distribution. We use the copula to construct a multivariate Meixner
process. Here the marginal density functions are given by
\begin{eqnarray}
\begin{array}{rl}
f_{M_i(t)}(x_i;m_i,a_i)=&\frac{(2cos(a_i/2))^{2m_i}}{2\pi\Gamma(2m_i)}
exp(a_ix_i)|\Gamma(m_i+\textrm{i}x_i)|^2,\\
&x_i\in(-\infty,+\infty),\quad i=1,2,\cdots,n
\end{array}\nonumber
\end{eqnarray}
The corresponding distribution is the measure of orthogonality of
the Meixner-Pollaczek polynomials (Koekoek and Swarttouw, 1998). The
Meixner process was introduced in Schoutens and Teugels (1998). In
Grigelions (1998), it is proposed for a model for risky assets and
an analogue of the famous Black and Scholes formula in mathematical
finance was established. The marginal characteristic functions of
$M_i(t)$, $i=1,2,\cdots,n$, are given by
\begin{eqnarray}
E\left[exp(\theta_iM_i(t))\right]=\left(\frac{cos(a_i/2)}{cosh((\theta_i-\textrm{i}a_i)/2)}\right)^{2m_it},
\quad i=1,2,\cdots,n.\nonumber
\end{eqnarray}
In according to the results in Schoutens and Teugels (1998) and
Schoutens(1999) and applying (11), its L\'{e}vy measure can be
calculated as:
\begin{eqnarray}
v(d\verb"x")&=&\partial_1\cdots\partial_nF|_{\xi_1=U_1(x_1),\cdots,\xi_n=U_n(x_n)}\prod\limits_{i=1}^n\frac{m_iexp(a_ix_i)}{x_isinh(\pi
x_i)}dx_i\nonumber\\
&=&\partial_1\cdots\partial_nF|_{\xi_1=U_1(x_1),\cdots,\xi_n=U_n(x_n)}\prod\limits_{i=1}^nm_i|\Gamma(1+\textrm{i}x_i)|^2\frac{exp(a_ix_i)}{\pi
x_i^2}dx_i.\nonumber
\end{eqnarray}
Also note that
\begin{eqnarray}
\left(\prod\limits_{i=1}^nx_i^2\right)v(d\verb"x")=\partial_1\cdots\partial_nF|_{\xi_1=U_1(x_1),\cdots,\xi_n=U_n(x_n)}
\prod\limits_{i=1}^nm_i|\Gamma(1+\textrm{i}x_i)|^2
\frac{exp(a_ix_i)}{\pi}dx_i.\nonumber
\end{eqnarray}
Being completely similar as in the above two examples, we can
orthogonalize the multivariate Teugels martingales for the
multivariate Meixner process by isometry.

\noindent{\textbf{Acknowledgment.}} The author thank a kind proposal
given by Professor David Nualart for the initial version of this
paper.

\end{document}